\documentclass[fleqn]{article}
    \usepackage{amsmath,amssymb,amsthm,url}
    \usepackage{graphicx}
    \usepackage{hyperref}
    \usepackage{amsrefs}

\title{Constructions of Generalized Sidon Sets}
\author {Greg Martin\thanks{University of British Columbia, \url{gerg@math.ubc.ca}}\
        \ and\
        {Kevin O'Bryant\thanks{University of California at San Diego, \url{kevin@member.ams.org}}}
        }
\date{\today}

 \newcommand{\MathReview}[1]{~\href{http://www.ams.org/mathscinet-getitem?mr=#1}{{\bf MR~#1}}}

 \newcommand{\ZZ}{{\mathbb Z}}

 \newcommand{\G}[2]{\ensuremath{\| #1 ^{\ast}\|_\infty}}

 \newcommand{\cK}{{\cal K}}

 \newcommand{\bigO}[1]{O\left(#1\right)}

 \newcommand{\floor}[1]{\left\lfloor #1 \right\rfloor}
 \newcommand{\ceiling}[1]{\big\lceil #1 \big\rceil}
 
 \newcommand{\field}[1]{{\mathbb F}_{#1}} 

 \newcommand{\Ruzsa}[3]{{\tt{Ruzsa}}(#1,#2,#3)}
 \newcommand{\Bose}[3]{{\tt{Bose}}(#1,#2,#3)}
 \newcommand{\Singer}[3]{{\tt{Singer}}(#1,#2,#3)}

 \newtheorem{thm}{Theorem}

 \newenvironment{proofof}[1]{\medskip\noindent{\em Proof of #1.}}{\qed\medskip}

\begin{document}
    \maketitle

\begin{abstract}
We give explicit constructions of sets $S$ with the property that for each integer $k$, there are at most $g$
solutions to $k=s_1+s_2, s_i\in S$; such sets are called Sidon sets if $g=2$ and generalized Sidon sets if $g\ge
3$. We extend to generalized Sidon sets the Sidon-set constructions of Singer, Bose, and Ruzsa. We also further
optimize Kolountzakis' idea of interleaving several copies of a Sidon set, extending the improvements of
Cilleruelo \& Ruzsa \& Trujillo, Jia, and Habsieger \& Plagne. The resulting constructions yield the largest
known generalized Sidon sets in virtually all cases.
\end{abstract}


\section{Sidon's Problem}
\label{section.Introduction}

In connection with his study of Fourier series, Simon Sidon~\cite{1932.Sidon} was led to ask how dense a set of
integers can be without containing any solutions to
    \begin{equation*}\label{SidonEq}
    s_1+s_2=s_3+s_4
    \end{equation*}
aside from the trivial solutions $\{s_1,s_2\}=\{s_3,s_4\}$. This, and certain generalizations, have come to be
known as {\em Sidon's Problem}.

Given a set $S\subseteq\ZZ$, we define the function $S\ast S$ by
    \begin{equation*}
        S\ast S (k) := \big| \{ (s_1,s_2) \colon s_i \in S, s_1+s_2=k\} \big|,
    \end{equation*}
which counts the number of ways to write $k$ as a sum of two elements of $S$. We also set
    \begin{equation*}
    \G{S}{2}
    := \|S\ast S\|_\infty = \max_{k\in\ZZ}\big| \{ (s_1,s_2) \colon s_i \in S, \, s_1+s_2=k \} \big|.
    \end{equation*}
Note that if the set $S$ is translated by $c$, then the function $S\ast S$ is translated by $2c$, and $\G{S}{2}$
is unaffected. Similarly, if the set $S$ is dilated by a factor of $c$, then $\G{S}{2}$ is unaffected.

If $\G{S}{2} \le 2$, then $S$ is called a Sidon set. Table~\ref{Sidon.Table} contains the optimally dense Sidon
sets with 10 or fewer elements. Erd\H{o}s \& Tur\`{a}n~\cite{1941.Erdos.Turan} showed that if
$S\subseteq[n]:=\{1,2,\dots,n\}$ is a Sidon set, then $|S|<n^{1/2}+\bigO{n^{1/4}}$, and
Singer~\cite{1938.Singer} gave a construction that yields a Sidon set in $[n]$ with $|S|>n^{1/2}-n^{5/16}$, for
sufficiently large $n$. Thus, the maximum density of a finite Sidon set is asymptotically known. The maximum
growth rate of $|S\cap [n]|$ for an infinite Sidon set $S$ remains enigmatic. We direct the reader
to~\cite{2004.OBryant} for a survey of Sidon's Problem.

\begin{table}[ht]\footnotesize \begin{center}\label{Sidon.Table}
\begin{tabular}{ccc}
$k$  & $\min\{a_k-a_1\}$ &             Witness              \\
\hline
2  &        1        &             \{0,1\}              \\
3  &        3        &            \{0,1,3\}             \\
4  &        6        &           \{0,1,4,6\}            \\
5  &       11        &          \{0,1,4,9,11\}          \\
  &                 &          \{0,2,7,8,11\}          \\
6  &       17        &       \{0,1,4,10,12,17\}         \\
  &                 &       \{0,1,4,10,15,17\}         \\
  &                 &       \{0,1,8,11,13,17\}         \\
  &                 &       \{0,1,8,12,14,17\}         \\
7  &       25        &      \{0,1,4,10,18,23,25\}       \\
  &                 &      \{0,1,7,11,20,23,25\}       \\
  &                 &      \{0,1,11,16,19,23,25\}      \\
  &                 &      \{0,2,3,10,16,21,25\}       \\
  &                 &      \{0,2,7,13,21,22,25\}       \\
8  &       34        &     \{0,1,4,9,15,22,32,34\}      \\
9  &       44        &   \{0,1,5,12,25,27,35,41,44\}    \\
10 &       55        &  \{0,1,6,10,23,26,34,41,53,55\}  \\
\end{tabular}
\caption{shortest Sidon sets, up to translation and reflection}
\end{center}\end{table}

The object of this paper is to give constructions of large finite sets $S$ satisfying the constraints
$S\subseteq[n]$ and $\G{S}{2}\le g$, that is, ``large'' in terms of $n$ and $g$. We extend the Sidon set
construction of Singer, as well as those of Bose~\cite{1942.Bose} and Ruzsa~\cite{1993.Ruzsa}, to allow
$\G{S}{2}\le g$ for arbitrary $g$. The essence of our extension is that although the union of 2 distinct Sidon
sets typically has large $\G{S}{2}$, the union of two of Singer's sets will have $\G{S}{2}\le 8$. We also
further optimize the idea of Kolountzakis~\cite{Kolountzakis} (refined in~\cite{Cilleruelo.Ruzsa.Trujillo} and
in~\cite{Habsieger.Plagne}) of controlling $\G{S}{2}$ by interleaving several copies of the {\em same} Sidon
set.

We warn the reader that the notation $\G{S}{2}$ is not in wide use. Most authors write ``$S$ is a $B_2[g]$
set'', sometimes meaning that $\G{S}{2}\le 2g$ and sometimes that $\G{S}{2} \le 2g+1$. Our notation is motivated
by the common practice of using the same symbol for a set and for its indicator function. With this convention,
    \[
    S\ast S(k) = \sum_{x\in \ZZ} S(x)S(k-x)
    \]
is the Fourier convolution of the function $S$ with itself, and counts representations as a sum of two elements
of $S$. We use the same notation when discussing subsets of $\ZZ_n$, the integers modulo $n$, and no ambiguity
arises.

\begin{figure}
    \begin{center}
    \begin{picture}(384,240)
        \put(376,120){$g$}
        \put(0,0){\includegraphics{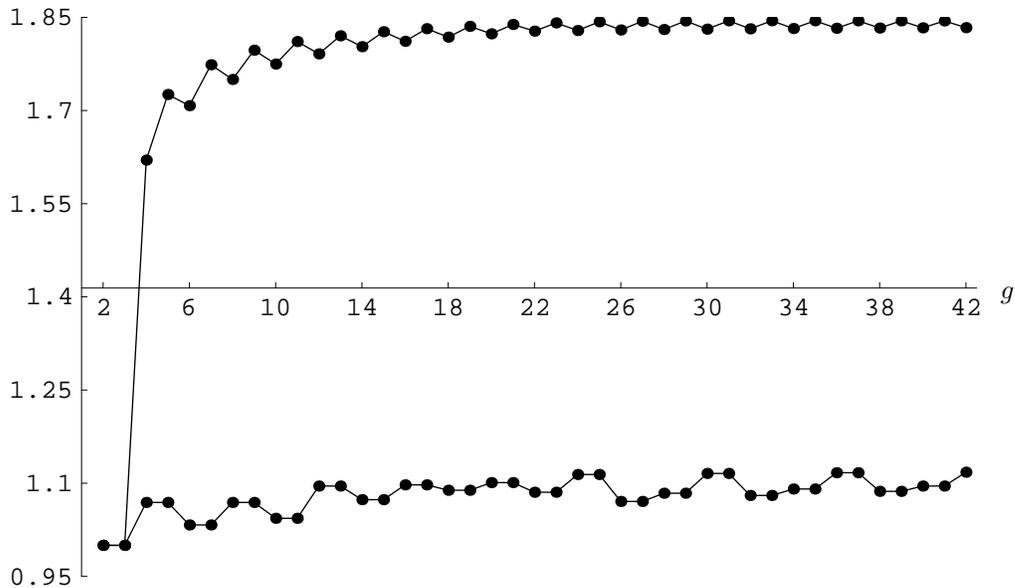}}
    \end{picture}
    \end{center}
    \caption{Upper and lower bounds on ${\sigma}(g)$.}
    \label{BothBoundsOnSigma.pic}
\end{figure}

Define
    \begin{equation}\label{Rdef}
    R(g,n):= \max_{S} \left\{ |S| \colon S\subseteq [n],\, \G{S}{2}\le g \right\}.
    \end{equation}
In words, $R(g,n)$ is the largest possible size of a subset of $[n]$ whose pairwise sums repeat at most $g$
times. We provide explicit lower bounds on $R(g,n)$ which are new for large values of $g$.
Figure~\ref{BothBoundsOnSigma.pic} shows the current upper and lower bounds on
    \[
    {\sigma}(g):=\liminf_{n\to\infty} \frac{R(g,n)}{\sqrt{\floor{g/2}n}}.
    \]
We comment that it may be possible to replace the $\liminf$ in the definition of $\sigma$ with a simple $\lim$,
but that this has not been proven and is not important for the purposes of this paper. The lower bounds on
$\sigma(g)$ all presented in this paper; some are originally found in~\cite{1938.Singer} ($g=2,3$),
\cite{Habsieger.Plagne} ($g=4$), and~\cite{Cilleruelo.Ruzsa.Trujillo} ($g=8,10$) but for other $g$ are new.
Other than the precise asymptotics for the $g=2$ and $g=3$ cases (which were found in 1944~\cite{1944.Erdos} and
1996~\cite{1996.Ruzsa}), the upper bounds indicated in Figure~\ref{BothBoundsOnSigma.pic} are due to
Green~\cite{2001.Green} when $g\le20$ is even; for all other values of $g$, the upper bounds are new and are the
subject of a work in progress by the authors~\cite{Martin.OBryant}.

Essential to proving these bounds on ${\sigma}(g)$ is the consideration of
    \begin{equation}\label{Cdef}
    C(g,n):= \max_{S} \left \{ |S| \colon S\subseteq \ZZ_n, \, \G{S}{2}\le g \right\}.
    \end{equation}
The function $C(g,n)$ gives the largest possible size of a subset of the integers modulo $n$ whose pairwise sums
(mod~$n$) repeat at most $g$ times. There is a sizable literature on $R(g,n)$, but little work has been done on
$C(g,n)$. There is a growing consensus among researchers on Sidon's Problem that substantial further progress on
the growth of $R(g,n)$ will require a better understanding of $C(g,n)$. Theorems~\ref{C.upperbound}
and~\ref{C.constructions} below give the state-of-the-art upper and lower bounds.

Tables~\ref{R.Table} and~\ref{C.Table} contain exact values for $R(g,n)$ and $C(g,n)$, respectively, for small
values of $g$ and $n$. These tables have been established by direct (essentially exhaustive) computation.
Specifically, Table~\ref{R.Table} records, for given values of $g$ and $k$, the smallest possible value of $\max
S$ given that $S\subseteq \ZZ^+$, $|S|=k$ and $\G{S}{2}\le g$; in other words, the entry corresponding to $k$
and $g$ is $\min\{n\colon R(g,n) \geq k\}$. For example, the fact that the $(k,g)=(8,2)$ entry equals 35 records
the fact that there exists an 8-element Sidon set of integers from $[35]$ but no 8-element Sidon set of integers
from $[34]$.

To show that $R(2,35)\geq 8$, for instance, it is only necessary to observe that
    \[
    S=\{1, 3, 13, 20, 26, 31, 34, 35\}
    \]
has 8 elements and $\G{S}{2}=2$. To show that $R(2,35)\leq 8$, however, seems to require an extensive search.

\begin{table}\footnotesize\begin{center}
$g$\vskip4pt $k$\quad
\begin{tabular}{|c||r|r|r|r|r|r|r|r|r|r|}
\hline
    &          2 &          3 &         4 &         5 &         6 & 7 & 8 & 9 &    10     &    11     \\ \hline\hline
 3  &          4 &            &           &           &           & & & & &           \\ \hline
 4  &          7 &          5 &           &           &           & & & & &           \\ \hline
 5  &         12 &          8 &         6 &           &           & & & & &           \\ \hline
 6  &         18 &         13 &         8 &         7 &           & & & & &           \\ \hline
 7  &         26 &         19 &        11 &         9 &         8 & & & & &           \\ \hline
 8  &         35 &         25 &        14 &        12 &        10 & 9 & & & &           \\ \hline
 9  &         45 &         35 &        18 &        15 &        12 & 11 & 10 &           &           &           \\ \hline
 10 &         56 &         46 &        22 &        19 &        14 & 13 & 12 &        11 &           &           \\ \hline
 11 &         73 &         58 &        27 &        24 &        17 & 15 & 14 &        13 &    12     &           \\ \hline
 12 &  $\leq 92$ &  $\leq 72$ &        31 &        29 &        20 & 18 & 16 &        15 &    14     &        13 \\ \hline
 13 &            &            &        37 &        34 &        24 & 21 & 18 &        17 &    16     &        15 \\ \hline
 14 &            &            &        44 &        40 &        28 & 26 & 21 &        19 &    18     &        17 \\ \hline
 15 &            &            & $\leq 52$ & $\leq 47$ &        32 & 29 & 24 &        22 &    20     &        19 \\ \hline
 16 &            &            &           &           &        36 & 34 & 27 &        24 &    22     &        21 \\ \hline
 17 &            &            &           &           & $\leq 42$ & $\leq 38$ &        30 &        28 &    24     &        23 \\ \hline
 18 &            &            &           &           &           & & 34 & 32 &    27     &        25 \\ \hline
 19 &            &            &           &           &           & & $\leq 38$ & $\leq 36$ &    30     &        28 \\ \hline
 20 &            &            &           &           &           & & & & 33     &        31 \\ \hline
 21 &            &            &           &           &           & & & & $\leq 37$ &        35 \\ \hline
 21 &            &            &           &           &           & & & & & $\leq 38$ \\ \hline
\end{tabular}
\caption{$\min\{n\colon R(g,n) \geq k\}$}\label{R.Table}
\end{center}\end{table}

\begin{table}
\footnotesize\begin{center} $g$\vskip4pt $k$\quad
\begin{tabular}{|c||r|r|r|r|r|r|r|r|r|r|}
\hline
    &         2 &  3 &  4 &  5 &  6 &  7 &  8 &  9 & 10 & 11 \\ \hline \hline
  3 &         6 &    &    &    &    &    &    &    &    &    \\ \hline
  4 &        12 &  7 &    &    &    &    &    &    &    &    \\ \hline
  5 &        21 & 11 &  8 &    &    &    &    &    &    &    \\ \hline
  6 &        31 & 19 & 11 &  9 &    &    &    &    &    &    \\ \hline
  7 &        48 & 29 & 14 & 13 & 10 &    &    &    &    &    \\ \hline
  8 &        57 & 43 & 22 & 17 & 12 & 11 &    &    &    &    \\ \hline
  9 &        73 & 57 & 28 & 19 & 16 & 13 & 12 &    &    &    \\ \hline
 10 &        91 &    & 36 & 28 & 19 & 17 & 14 & 13 &    &    \\ \hline
 11 &           &    &    & 35 & 22 & 21 & 18 & 15 & 14 &    \\ \hline
 12 &           &    &    &    & 30 & 23 & 21 & 19 & 16 & 15 \\ \hline
 13 &           &    &    &    &    & 31 & 24 & 22 & 19 & 17 \\ \hline
 14 &           &    &    &    &    &    & 28 & 25 &    & 20 \\ \hline
\end{tabular}
\caption{$\min\{n\colon C(g,n) \geq k\}$}\label{C.Table}
\end{center}\end{table}

In the next section, we state our upper bounds on $C(g,n)$, lower bounds on $R(g,n)$ and $C(g,n)$, and
constructions that demonstrate our lower bounds. In Section~\ref{section.Proofs} we prove the bounds claimed in
Section~\ref{section.Theorems}. Since the value of this work is primarily as a synthesis and extension of ideas
from a variety of other works, we have endeavored to make this paper self-contained. We conclude in the final
section by listing some questions that we would like, but have been unable, to answer.

\section{Theorems and Constructions}
\label{section.Theorems}

\subsection{Theorems}

\begin{thm} \label{C.upperbound}
    \renewcommand{\theenumi}{\roman{enumi}}
    \begin{enumerate}
    \item[(i)]   $\binom{C(2,n)}{2} \le \floor{\frac n2}$, and in particular $C(2,n)\le \sqrt{n}+1$;
    \item[(ii)]   $C(3,n) \le \sqrt{n+9/2}+3$;
    \item[(iii)]   $C(4,n) \le \sqrt{3n} + 7/6$;
    \item[(iv)]   $C(g,n) \le \sqrt{gn}$ for even $g$;
    \item[(v)]   $C(g,n) \le \sqrt{1-\tfrac1g}\sqrt{gn}+1$, for odd $g$.
    \end{enumerate}
\end{thm}

\begin{thm}
    \label{C.constructions}
    \renewcommand{\theenumi}{\roman{enumi}}
    Let $q$ be a prime power, and let $k, g,f,x,y$ be positive integers with $k<q$.
    \begin{enumerate}
    \item[(i)] If $p$ is a prime, then $C(2k^2,p^2-p) \geq k (p-1)$;
    \item[(ii)] $C(2k^2,q^2-1) \geq k q$;
    \item[(iii)] $C(2k^2,q^2+q+1) \geq kq+1$;
    \item[(iv)] If $\gcd(x,y)=1$, then $C(gf,xy)\geq C(g,x)C(f,y)$;
    \item[(v)] $R(gf,xy) \geq R(gf,xy+1-\ceiling{\tfrac{y}{C(f,y)}}) \geq R(g,x)C(f,y)$;
    \item[(vi)] $R(g,3g-\floor{g/3}+1) \geq g+2\floor{g/3}+\floor{g/6}$.
    \end{enumerate}
\end{thm}

\begin{thm}
\label{ExplicitValues}
$$
\begin{array}{r@{{}\ge{}}l@{{}>{}}l}
    \sigma (4)  & \sqrt{8/7}                   & 1.069, \\ \vspace{1mm}
    \sigma (6)  & \sqrt{16/15}                 & 1.032, \\ \vspace{1mm}
    \sigma (8)  & \sqrt{8/7}                   & 1.069, \\ \vspace{1mm}
    \sigma (10) & \sqrt{49/45}                 & 1.043, \\ \vspace{1mm}
    \sigma (12) & \sqrt{6/5}                   & 1.095,
\end{array}
\qquad
\begin{array}{r@{{}\ge{}}l@{{}>{}}l}
    \sigma (14) & \sqrt{121/105}               & 1.073, \\ \vspace{1mm}
    \sigma (16) & \sqrt{289/240}               & 1.097, \\ \vspace{1mm}
    \sigma (18) & \sqrt{32/27}                 & 1.088, \\ \vspace{1mm}
    \sigma (20) & \sqrt{40/33}                 & 1.100, \\ \vspace{1mm}
    \sigma (22) & \sqrt{324/275}               & 1.085,
\end{array}
$$
\end{thm}

\begin{thm}
\label{R.Lower.Bounds.thm} For $g\ge1$,
   \[
   \sigma(2g+1) \geq \sigma (2g) \ge \frac{g+2\floor{g/3}+\floor{g/6}}{\sqrt{3g^2-g\floor{g/3}+g}}\,.
   \]
In particular,
    \[
    \liminf_{g\to\infty} \sigma(g) \ge \frac{11}{\sqrt{96}}.
    \]
\end{thm}

We note that Martin \& O'Bryant have shown~\cite{Martin.OBryant} that $\limsup_{g\to\infty} \sigma(g) < 1.8391$,
whereas $11/\sqrt{96}>1.1226$. These lower bounds on $\sigma$, together with the strongest known upper bounds,
are plotted for $2\leq g \leq 42$ in Figure~\ref{BothBoundsOnSigma.pic}.

\subsection{Constructions}

Theorem~\ref{C.constructions} rests on the constructions given in the following four subsubsections. We denote
the finite field with $q$ elements by $\field{q}$, and its multiplicative group by $\field{q}^\times$.

\subsubsection{Ruzsa's Construction}
Let $\theta$ be a generator of the multiplicative group modulo the prime $p$. For $1\leq i < p$, let $a_{t,i}$
be the congruence class modulo $p^2-p$ defined by
    \begin{equation*}
    a_{t,i} \equiv t \pmod{p-1} \quad \text{ and } \quad
    a_{t,i} \equiv i \theta^t \pmod{p}.
    \end{equation*}
Define the set
    \[
    \Ruzsa{p}{\theta}{k} := \{a_{t,k} \colon 1\le t < p \} \subseteq \ZZ_{p^2-p}.
    \]
Ruzsa~\cite{1993.Ruzsa} showed that $\Ruzsa{p}{\theta}{1}$ is a Sidon set. We show that if ${\cal K}$ is any
subset of $[p-1]$, then
    \[
    \Ruzsa{p}{\theta}{\cK} := \bigcup_{k\in \cK} \Ruzsa{p}{\theta}{k}
    \]
is a subset of $\ZZ_{p^2-p}$ with cardinality $|{\cK}|(p-1)$ and
    \[\G{\Ruzsa{p}{\theta}{\cK}}{2} \le {2|{\cK}|^2}.\]

For example, $\Ruzsa{11}{2}{\{1,2\}}$ is
    \[
                             \{7, 39, 58, 63, 65, 86, 92, 100, 101, 104 \}
                                \cup
                             \{28, 47, 52, 54, 75, 81, 89, 90, 93, 106\},
    \]
and one may directly verify that $\G{\Ruzsa{11}{2}{\{1,2\}}}{2} =8$.

\subsubsection{Bose's Construction}

Let $q$ be any prime power, $\theta$ a generator of $\field{q^2}$, $k \in \field{q}$, and define the set
    \[
    \Bose{q}{\theta}{k} := \{a \in [q^2-1] \colon \theta^a-k\theta\in \field{q}\}.
    \]
Bose~\cite{1942.Bose} showed that for $k\not=0$, $\Bose{q}{\theta}{k}$ is Sidon set. We show that if ${\cal K}$
is any subset of $\field{q}\setminus\{0\}$, then
    \[
    \Bose{q}{\theta}{\cK} := \bigcup_{k\in \cK} \Bose{q}{\theta}{k}
    \]
is a subset of $\ZZ_{q^2-1}$, has $|\cK|q$ elements, and
    \[
    \G{\Bose{q}{\theta}{\cK}}{2} \le 2 | \cK |^2.
    \]

For example, $\Bose{11}{x \bmod{(11,x^2+3x+6)}}{\{1,2\}}$ is
    \[
                           \{1, 30, 38, 55, 56, 65, 69, 71, 76, 99, 118\}
                                \cup
                           \{18, 26, 43, 44, 53, 57, 59, 64, 87, 106, 109\}.
    \]

\subsubsection{Singer's Construction}

Sidon sets arose incidentally in Singer's work~\cite{1938.Singer} on finite projective geometry. While Singer's
construction gives a slightly thicker Sidon set than Bose's (which is slightly thicker than Ruzsa's), the
construction is more complicated---even after the simplification of~\cite{1962.Bose.Chowla}.

Let $q$ be any prime power, and let $\theta$ be a generator of the multiplicative group of $\field{q^3}$. For
each $k_1,k_2 \in\field{q}$ define the set
    \[
    T(\langle k_1,k_2 \rangle)
    := \{0\} \cup \{ a \in [q^3-1] \colon \theta^a- k_2 \theta^2-k_1\theta \in \field{q}\}.
    \]
Then define
    \[
    \Singer{q}{\theta}{\langle k_1,k_2 \rangle}
    \]
to be the congruence classes modulo $q^2+q+1$ that intersect $T(\langle k_1,k_2 \rangle)$. Singer proved that
for $k_2=0,k_1\not=0$, $\Singer{q}{\theta}{\langle k_1,k_2 \rangle}$ is a Sidon set. We show that if
$\cK\subseteq \field{q}\times\field{q}$ does not contain two pairs with one an $\field{q}$-multiple of the
other, then
    \[
    \Singer{q}{\theta}{\cK} := \bigcup_{\vec{k}\in \cK} \Singer{q}{\theta}{\vec k}
    \]
is a subset of $\ZZ_{q^2+q+1}$ with $|\cK| q +1$ elements and
    \[
    \G{\Singer{q}{\theta}{\cK}}{2} \le 2 | \cK |^2.
    \]

For example, $\Singer{11}{x\bmod{(11,x^3 + x^2 + 6x + 4)}}{\{\langle 1,1\rangle,\langle 1,2\rangle\}}$ is
    \[
    \{0, 9, 57, 59, 63, 81, 86, 97, 100, 112, 125, 132\}
    \cup
    \{3, 15, 28, 35, 36, 45, 93, 95, 99, 117, 122\}.
    \]

\subsubsection{The Cilleruelo \& Ruzsa \& Trujillo Construction}

Kolountzakis observed that if $S$ is a Sidon set, and $S+1:=\{s+1\colon s\in S\}$, then
$\G{(S\cup(S+1))}{2}\le4$. This idea of interleaving several copies of the same Sidon set was extended
incorrectly by Jia (but fixed by Lindstr\"{o}m), and then correctly by Cilleruelo \& Ruzsa \& Trujillo, Habsieger \&
Plagne, and Cilleruelo (to $h>2$).

Let $S\subseteq \ZZ_{x}$ and $M\subseteq \ZZ_{y}$ have $\G{S}{2}\le g$ and $\G{M}{2}\le f$. Let $S' \subseteq
[x]$ and $M' \subseteq [y]$ be corresponding sets of integers, i.e., $S=\{s \mod{x} \colon s\in S'\}$. Now, let
    \[
    M'+yS' := \{ m+ys \colon m\in M', s\in S' \}\subseteq \ZZ.
    \]
The set
    \[
    M+yS := \{t \bmod{xy} \colon t\in M'+yS'\}\subseteq \ZZ_{xy}
    \]
satisfies $\G{(M+yS)}{2} \le gf$.

\section{Proofs}
\label{section.Proofs}

If $S$ is a set of integers (or congruence classes), we use $S(x)$ to denote the corresponding indicator
function. Also, we use the standard notations for convolution and correlation of two real-valued functions:
    \[
    S\ast T (x) = \sum_y S(y)T(x-y) \quad \text{ and } \quad  S\circ T(x) = \sum_y S(y)T(x+y).
    \]
For sets $S,T$ of integers, $S\ast T(x)$ is the number of ways to write $x$ as a sum $s+t$ with $s\in S$ and
$t\in T$. Likewise, $S\circ T(x)$, is the number of ways to write $x$ as a difference $t-s$.

\subsection{Theorem~\ref{C.upperbound}}

Part (i) is just the combination of the pigeonhole principle and the fact (which we prove below) that if
$\|S\ast S\|_\infty \le 2$, then for $k\not=0$, $S\circ S(k)\le 1$. Part (ii) follows from the observation that
if $\|S\ast S\|_\infty \le 3$, then $S\circ S(k)\le 2$ for $k\not=0$, and in fact $S\circ S(k)\le 1$ for almost
all $k$. Part (iii) follows an idea of Cilleruelo: if $\|S\ast S\|_\infty \le 4$, then $S\circ S$ is small on
average. For $g>4$, the theorem is a straightforward consequence of the pigeonhole principle. We consider part
(iii) to be the interesting contribution.

\proofof{(i)} We show that $\binom{C(2,n)}{2}\leq \floor{n/2}$, whence $C(2,n)<\sqrt{n}+1$. Let $S\subseteq[n]$
have $\G{S}{2}\le 2$. If $\{s_1,s_2\}$, $\{s_3,s_4\}$ are distinct pairs of distinct elements of $S$, and
    \begin{equation}\label{equaldifferences}
    s_1-s_2 \equiv s_3-s_4 \pmod{n},
    \end{equation}
then $s_4+s_1 \equiv s_1+s_4 \equiv s_3+s_2 \equiv s_3+s_2$, contradicting the supposition that $\G{S}{2}\le2$.
Therefore, the map $\{s_1,s_2\}\mapsto \{\pm(s_1-s_2)\}$ is 1-1 on pairs of distinct elements of $S$, and the
image is contained in $\{\{\pm1\}, \{\pm2\}, \dots, \{\pm\floor{n/2}\}\}$. Thus, $\tbinom{|S|}{2}
\le\floor{n/2}$.

This bound is actually achieved for $n=p^2+p+1$ when $p$ is prime (see Theorem~\ref{C.constructions}(iii)).

\proofof{(ii)} Now suppose that $\G{S}{2} = 3$, and consider the pairs of distinct elements of $S$. Any solution
to \eqref{equaldifferences} must have $\{s_1,s_2\}\cap\{s_3,s_4\}\not= \emptyset$ since $\G{S}{2}<4$, and each
of the $|S|$ possible intersections can occur only once. Therefore, after deleting one pair for each element of
$S$, we get a set of $\tbinom{|S|}2-|S|$ pairs which is mapped 1-1 by $\{s_1,s_2\}\mapsto \{\pm(s_1-s_2)\}$ into
$\{\{\pm1\}, \{\pm2\}, \dots, \{\pm\floor{n/2}\}\}$. This proves Theorem~\ref{C.upperbound} for $g=3$.

\proofof{(iii)} Now suppose that $\G{S}{2}=4$, where $S\subseteq\ZZ_{n}$. The obvious map from
    \begin{align*}
    X &:= \left\{ \Big( (s_1,s_2),(s_3,s_4)\Big) \colon s_1-s_2 \equiv s_3-s_4, s_1
                \not\in\{s_2,s_3\}\right\} \\
\intertext{to}
    Y &:= \left\{ \Big( (s_1,s_4),(s_3,s_2)\Big) \colon s_1+s_4 \equiv s_3+s_2,
                \{s_1,s_4\} \not\equiv \{s_2,s_3\} \right\}
    \end{align*}
is easily seen to be 1-1 (but not necessarily onto): $|X|\le |Y|$. We have
    \begin{multline*}
    |X| =  \sum_{\substack{k\not=0 \\ k \in \ZZ_{n}}} \left(S\circ S(k)^2-S\circ S(k)\right)
        \ge \frac{1}{n-1} \bigg( \sum_{\substack{k\not=0 \\ k \in \ZZ_{n}}} S\circ S(k) \bigg)^2
                -\sum_{\substack{k\not=0 \\ k \in \ZZ_{n}}} S\circ S(k) \\
        = \frac{(|S|^2-|S|)^2}{n-1}  - |S|^2+|S|
    \end{multline*}
    \begin{equation*}
    |Y| =  \left|(S\ast S)^{-1}(3)\right|\, 4 + \left|(S\ast S)^{-1}(4)\right| \,8
        \le4|S| + 8 \frac{|S|^2-|S|}{4}
        = 2|S|^2+2|S|
    \end{equation*}
Comparing the lower bound on $|X|$ with the upper bound on $|Y|$ yields $|S| \le \sqrt{3n}+7/6$.

\proofof{(iv) and (v)} There are $|S|^2$ pairs of elements from $S\subseteq\ZZ_{n}$, and there are just $n$
possible values for the sum of two elements. If $\G{S}{2}\le g$ then each possible value is realized at most $g$
times. Thus $|S|^2 \leq {gn}$. The only way a sum can occur an odd number of times is if it is twice an element
of $S$, so for odd $g$, $|S|^2 \leq (g-1)n+|S|$.

\subsection{Theorem~\ref{C.constructions}}

The first three parts of Theorem~\ref{C.constructions} are all proved in a similar manner, which we outline
here. For disjoint sets $S_1, \dots S_k$, with $S=\cup S_i$, we have
    \[
    S\ast S = (S_1+\dots+S_k) \ast (S_1+\dots +S_k) = \sum_{i,j=1}^k S_i \ast S_j
    \]
and since $S_i \ast S_j$ is nonnegative,
    \[
    \|S\ast S\|_\infty \le \sum_{i,j=1}^k \|S_i \ast S_j\|_\infty \le k^2 \max_{1\le i,j\le k} \|S_i \ast
    S_j\|_\infty.
    \]

To prove part (i), we need to show that the sets $\Ruzsa{p}{\theta}{i}$ ($1\le i <p$) are disjoint (hence
$\Ruzsa{p}{\theta}{\cK}$ has cardinality $|\cK|(p-1)$), and that
    \[
    \|\Ruzsa{p}{\theta}{i} \ast \Ruzsa{p}{\theta}{j} \|_\infty \le 2.
    \]
Specifically, we use unique factorization in $\field{p}[x]$ to show that there are not 3 distinct pairs
    \[(a_{r_m,i},a_{v_m,j}) \in \Ruzsa{p}{\theta}{i} \times \Ruzsa{p}{\theta}{j}\]
with the same sum.

The proofs of parts (ii) and (iii) follow the same outline, but use unique factorization in $\field{q}[x]$ and
$\field{q^2}[x]$, respectively.

\proofof{(i)} For the entirety of the proof, we work with fixed $p$ and $\theta$. It is therefore convenient to
introduce the notation $R_k = \Ruzsa{p}{\theta}{k}$. We need to show that $R_i \cap R_j =\emptyset $ for $1\le i
< j <p$, and that $\|R_i \ast R_j\|_\infty \le 2$ (including the possibility $i=j$).

Suppose that $a_{m_1,i}=a_{m_2,j} \in R_i \cap R_j$, with $m_1,m_2 \in [1,p)$. We have $m_1 \equiv a_{m_1,i}
=a_{m_2,j}\equiv m_2 \pmod{p-1}$, so $m_1=m_2$. Reducing the equation $a_{m_1,i}=a_{m_2,j}$ modulo $p$, we find
$i\theta^{m_1}\equiv j\theta^{m_2}=j\theta^{m_1} \pmod{p}$, so $i=j$. Thus for $i\not=j$, the sets $R_i,R_j$ are
disjoint.

Now suppose, by way of contradiction, that there are three pairs $(a_{r_m,i},a_{v_m,j})\in R_i \times R_j$
satisfying $a_{r_m,i}+a_{v_m,j}\equiv k \pmod{p^2-p}$. Each pair gives rise to a factorization modulo $p$ of
    \[
    x^2 - k x + ij \theta^k\equiv (x-a_{r_m,i})(x-a_{v_m,j}) \pmod{p} .
    \]
Factorization modulo $p$ is unique, so it must be that two of the three pairs are congruent modulo $p$, say
    \begin{equation}\label{modp:eq}
    a_{r_1,i}\equiv a_{r_2,i} \pmod{p}.
    \end{equation}
In this case, $i\theta^{r_1} \equiv a_{r_1,i}\equiv a_{r_2,i} \equiv i \theta^{r_2} \pmod{p}$. Since $\theta$
has multiplicative order $p-1$, this tells us that $r_1 \equiv r_2 \pmod{p-1}$. Since $a_{r_m,i} \equiv r_m
\pmod{p-1}$ by definition, we have
    \begin{equation}\label{modp-1:eq}
    a_{r_1,i} \equiv a_{r_2,i} \pmod{p-1}.
    \end{equation}
Equations \eqref{modp:eq} and \eqref{modp-1:eq}, together with
    \[
    a_{r_1,i}+a_{v_1,j} \equiv k \equiv a_{r_2,i}+a_{v_2,j} \pmod{p^2-p}
    \]
imply that the first two pairs are identical, and so there are {\em not} three such pairs. Thus, for each $k\in
\ZZ_{n}$, we have shown that $R_i \ast R_j(k) \le 2$.

\proofof{(ii)} For $k\in \field{q}$, let $B_k=\Bose{q}{\theta}{k}$. We need to show that $|B_i|=q$, that $B_i
\cap B_j =\emptyset$ for distinct $i,j\in \field{q}\setminus\{0\}$, and that $\|B_i\ast B_j\|_\infty \le 2$
(including the possibility that $i=j$).

Since $\{\theta,1\}$ is a basis for $\field{q^2}$ over $\field{q}$, we can for each $s'\in [q^2-1]$ write
$\theta^{s'}$ as a linear combination of $\theta$ and $1$. We define $s$ (unprimed) to be the coefficient of 1,
i.e.,
    \[
    \theta^{s'}= i \theta+ s
    \]
for some $i$. In this proof, primed variables are integers between 1 and $q^2-1$, and unprimed variables are
elements of $\field{q}$. Note also that $a'=b'$ implies $a=b$, whereas $a=b$ does not imply $a'=b'$.

Since $\theta$ generates the multiplicative group, for $i\not=0$ each $s\in\field{q}$ has a corresponding $s'$,
so that $|B_i|=q$. Moreover, we know that $i \theta+s_1=j\theta+s_2$ implies that $i=j$ and $s_1=s_2$. In
particular, if $i\not=j$, then $B_i \cap B_j=\emptyset$. Thus $|\Bose{q}{\theta}{\cK}|= |\cK| q$.

We now fix $i$ and $j$ in $\field{p}\setminus\{0\}$ (not necessarily distinct), and show that $B_i\ast B_j(k)\le
2$ for $k\in\ZZ_{q^2-1}$. Define $c_1,c_2 \in \field{p}$ by $(ij)^{-1} \theta^{k'}-\theta^2=c_1\theta+c_2$, and
consider pairs $(r',v') \in B_i \times B_j$ with $r'+v' \equiv k' \pmod{q^2-1}$. We have
 \begin{multline*}
        c_1\theta+c_2 = (ij)^{-1} \theta^{k'}-\theta^2 =
        (ij)^{-1} \theta^{r'+v'}-\theta^2 =
        (ij)^{-1} \theta^{r'}\theta^{v'}-\theta^2 = \\
        (ij)^{-1} (i \theta+r)(j \theta+v)-\theta^2 =
        (i^{-1}r+j^{-1}v)\theta+i^{-1}r j^{-1}v.
 \end{multline*}
This means that $(a,b)=(i^{-1}r, j^{-1}v)$ is a solution to $x^2-c_1 x+c_2=(x-a)(x-b)$. By unique factorization
over finite fields, there are at most two such pairs. Thus, $B_i \ast B_j(k)\le 2$ and so $\|B_i \ast B_j
\|_\infty \le 2$.

\proofof{(iii)} We first note that $\theta^a$ and $\theta^b$ (for any integers $a,b$) are linearly dependent
over $\field{q}$ if and only if their ratio is in $\field{q}$. Since $\field{q}^\times$ is a subgroup of
$\field{q^3}^\times$, we see that $\field{q}=\{\theta^{x(q^2+q+1)} \colon x\in\ZZ\}$. Thus, we have the
following linear dependence criterion: $\theta^a$ and $\theta^b$ are linearly dependent if and only if $a\equiv
b \pmod{q^2+q+1}$.

Since $\{\theta^2,\theta,1\}$ is a basis for $\field{q^3}$ over $\field{q}$, we can for each $s'\in [q^3-1]$
write $\theta^{s'}$ as a linear combination of $\theta^2$, $\theta$ and $1$. We define $s$ (unprimed) to be the
coefficient of 1, i.e.,
    \[
    \theta^{s'}= k_2\theta^2+k_1 \theta+ s
    \]
for some $k_1,k_2$. In this proof, primed variables are integers between 1 and $q^3-1$, and unprimed variables
are elements of $\field{q}$. Note, as above, that $a'=b'$ implies $a=b$, whereas $a=b$ does not imply $a'=b'$.
We also define $\bar s$ to be the congruence class of the integer $s'$ modulo $q^2+q+1$.

For $\vec{k}=\langle k_1,k_2 \rangle \in \field{q}^2$ define
    \[
    T(\vec{k}) :=  \{ s' \in [q^3-1] : \theta^{s'} = s+k_1\theta+k_2\theta^2,\quad  s\in\field{q}\}
    \]
which also reiterates the connection between primed variables (such as $s'\in [q^3-1]$) and unprimed variables
(such as $s\in \field{q}$). Define $S(\vec{k})$ to be the set of congruence classes modulo $q^2+q+1$ that
intersect $T(\vec{k})$; as noted above, we denote the congruence class $s'\bmod{q^2+q+1}$ as $\bar s$. Let
$\cK=\{\vec{k}_1,\vec{k}_2,\dots\} \subseteq \field{q}\times\field{q}$ be a set that does not contain two pairs
with one being a multiple of the other. Let $S_1:=\{0\}\cup S(\vec{k}_1)$, and for $i>1$ let
$S_i:=S(\vec{k}_i)$.

We need to show that $|S_1|=q+1$, for $i>1$ that $|S_i|=q$, and for distinct $i$ and $j$, the sets $S_i$ and
$S_j$ are disjoint. This will imply that
    \[
    \Singer{q}{\theta}{\cK} = \bigcup_{i=1}^{|\cK|} S_i
    \]
has cardinality $|\cK| q +1$. All of these are immediate consequences of the fact that each element of
$\field{q^3}$ has a unique representation as an $\field{q}$-linear combination of $\theta^2$, $\theta$, and 1.

We will show that for any $i,j$ (not necessarily distinct) there are not three pairs $(\bar{r}_m ,\bar{v}_m) \in
S_i \times S_j$ with the same sum modulo $q^2+q+1$.

Suppose that $\vec{k}_i=\langle k_1,k_2\rangle$ and $\vec{k}_j=\langle \ell_1,\ell_2\rangle$. Set $K(r,z):=r+k_1
z + k_2 z^2$ and $L(v,z)=v+\ell_1 z + \ell_2 z^2$. Since
    \[
    \bar r_1 +\bar v_1 =\bar r_2  + \bar v_2  = \bar r_3 +\bar v_3
    \]
there are constants $c_2,c_3\in \field{q}$ such that $\theta^{r_1'+v_1'}=c_2\theta^{r_2'+v_2'} = c_3
\theta^{r_3'+v_3'}$, and since $\theta^{r'+v'}=\theta^{r'}\theta^{v'}=K(r,\theta)L(v,\theta)$, the polynomials
    \begin{align*}
    f_2(z)  &:= c_2 K(r_2,z)L(v_2,z)-K(r_1,z)L(v_1,z) \\
    f_3(z)  &:= c_3 K(r_3,z)L(v_3,z)-K(r_1,z)L(v_1,z)
    \end{align*}
both have $\theta$ as a root (we are assuming for the moment that none of $\bar v_m,\bar r_m$ are $\bar0$).

If $c_2=1$, then $f_2(z)$ is a quadratic with the cubic $\theta$ as a root: consequently $f_2(z)=0$ identically.
This gives three equations in the unknowns $r_1$, $v_1$, $r_2$, $v_2$, $k_1$, $k_2$, $\ell_1$, $\ell_2$. These
equations with the assumption that $\langle k_1,k_2\rangle$ is not a multiple of $\langle \ell_1,\ell_2\rangle$,
imply that $r_1=r_2$ and $v_1=v_2$. Thus $\theta^{r_1'}=\theta^{r_2'}$, and so $r_1'=r_2'$, and so
$(r_1',v_1')=(r_2',v_2')$, contrary to our assumption of distinctness. Similarly $c_3\not=1$ and $c_2\not=c_3$.

Now
    \[
    g(z):= (c_3-1)f_2(z)-(c_2-1)f_3(z)
    \]
is a quadratic with $\theta$ as a root. Setting its coefficients equal to 0 gives 3 equations:
    \begin{align*}
    0 &= {c_2}\,\left( {r_1}\,{v_1} -{r_2}\,{v_2} \right)+
            {c_3}\,\left({r_3}\,{v_3}- {r_1}\,{v_1}   \right) +
            {c_2}\,{c_3}\left( \,{r_2}\,{v_2} - {r_3}\,{v_3} \right)
    \\
    0 &= {c_2}\,\big({\ell_1}\,({r_1} - {r_2}) + {k_1}\,({v_1} -{v_2})\big) +
            {c_3}\,\big( {\ell_1}(r_3-r_1)+k_1(v_3-v_1)\big)  \\
    &\hspace{6cm}
            +{c_2}\,{c_3}\,\big( {\ell_1}\,\left( {r_2} - {r_3} \big)  +
                    {k_1}\,\left( {v_2} - {v_3} \right)\right)
    \\
    0 &= {c_2}\,\big({\ell_2}\,\left( r_1-r_2 \right)  + {k_2}\,\left( v_1-v_2 \right)\big) +
            {c_3}\,\big({\ell_2}\,\left( r_3-{r_1}  \right)  +
                            {k_2}\,\left( v_3-{v_1} \right)\big) \\
    &\hspace{6cm}
            + {c_2}\,{c_3}\big({\ell_2}\,\left( {r_2} - {r_3} \right)  +
                                {k_2}\,\left( {v_2} - {v_3} \right)\big)
    \end{align*}
When combined with our knowledge that $c_2,c_3$ are not 0, 1, or equal, and $\langle k_1,k_2\rangle$ not a
multiple of $\langle \ell_1,\ell_2\rangle$, this implies that the pairs $(\bar r_m, \bar v_m)$ are not distinct.

Now suppose that $\bar r_1=0$, $\bar v_1 \not =0$, and set
    \begin{align*}
    f_2(z)  &:= c_2 K(r_2,z)L(v_2,z)-L(v_1,z) \\
    f_3(z)  &:= c_3 K(r_3,z)L(v_3,z)-L(v_1,z).
    \end{align*}
We have $f_2(\theta)=f_3(\theta)=0$, and in particular
    \[
    g(z):=c_3 f_2(z)-c_2 f_3(z)
    \]
is a quadratic with $\theta$ as a root. Setting the coefficients of $g(z)$ equal to 0 yields equations which, as
before, with our assumptions about $c_2,c_3,k_1,k_2,\ell_1,\ell_2$, imply that the three pairs $(\bar r_m, \bar
v_m)$ are not distinct. The case $\bar r_1=\bar v_1=0$ is handled similarly. The case $\bar r_1=\bar v_2=0$ is
eliminated for distinct $i,j$ by the disjointness of $S_i$ and $S_j$, and for $i=j$ by the distinctness
assumption on the three pairs.

Thus there are not such $(\bar r_m, \bar v_m)$ ($1\le m \le 3$), whether none of these six variables are 0, one
of them is 0, or two of them are 0.

\proofof{(iv)} Consider $m_i,n_i\in M'$ and $s_i,t_i\in S'$ with
    \begin{equation}
        \label{mainsupposition}
    (m_1+ys_1)+(n_1+yt_1)\equiv
    \dots\equiv(m_{gh+1}+ys_{gf+1})+ (n_{gf+1}+yt_{gf+1}) \pmod{xy}.
    \end{equation}
We need to show that $m_i=m_j$, $s_i=s_j$, $n_i=n_j$, and $t_i=t_j$, for some distinct $i,j$. Reducing
Eq.~\eqref{mainsupposition} modulo $y$, we see that $m_1+n_1\equiv m_2+n_2\equiv \dots\equiv
m_{gf+1}+n_{gf+1}\pmod{y}$. Since $\G{M}{2}\le f$, we can reorder the $m_i,n_i,s_i,t_i$ so that
$m_1=m_2=\dots=m_{g+1}$ and $n_1=n_2=\dots=n_{g+1}$. Reducing Eq.~\eqref{mainsupposition} modulo $x$ we arrive
at
    \[
    y s_1+y t_1\equiv y s_2+ yt_2\equiv \dots\equiv y s_{g+1}+y t_{g+1} \pmod{x}
    \]
whence, since $\gcd(x,y)=1$,
    \[
    s_1+t_1\equiv s_2+t_2\equiv \dots \equiv s_{g+1}+t_{g+1} \pmod{x}.
    \]
The $s_i\bmod{x}$ and $t_i\bmod{x}$ are from $S$, and $\G{S}{2}\le g$, so that for some distinct $i,j$,
$s_i=s_j$ and $t_i=t_j$.

\proofof{(v)} Let $M \subseteq \ZZ_{y}$ have cardinality $C(f,y)$ and $\G{M}{2}\le f$. Set $M'=\{ m\in [y]
\colon m\bmod{y} \in M\}$. Let $S'\subseteq[0,r)$ have cardinality $R(g,r)$ and $\G{(S')}{2}\le g$. Set (with
$x>2r$) $S := \{ s\bmod{x} \colon s \in S'\}\subseteq \ZZ_{x}$. By the construction in part (iv) of this theorem
$M+yS \subseteq \ZZ_{xy}$ has
    \[
    \G{(M+yS)}{2} \le gf.
    \]
Since $M'+yS'\subseteq [y+yr]$ and $M'+yS' \equiv M+yS \pmod{xy}$, if $xy>2(y+yr)$ then
$\G{(M'+yS')}{2}=\G{(M+yS)}{2} \le gf$.

We can shift $M$ modulo $y$ without affecting $|M|$ or $\G{M}{2}$. Since there clearly must be two consecutive
elements of $M$ with difference at least $\ceiling{y/C(f,y)}$, we may assume that
    $M'\subseteq[y-\ceiling{y/C(f,y)}+1].$
Thus,
    \[
    M'+yS' \subseteq [y-\ceiling{y/C(f,y)}+1+y(r-1)] =[yr+1-\ceiling{y/C(f,y)}]
    \]
and
    \[
    |M'+yS'| = |M| \, |S'| = C(f,y) R(g,r).
    \]
This proves part (v).

The reader might feel that the part of the argument concerning the largest gap in $M$ is more trouble than it is
worth. We include this for two reasons. First, Erd\H{o}s~\cite{1994.Guy}*{Problem C9} offered \$500 for an
answer to the question, ``Is $R(2,n)=\sqrt{n}+\bigO{1}$?'' This question would be answered in the negative if
one could show, for example, that $\Bose{p}{\theta}{1}$ contains a gap that is not $\bigO{p}$, as seems likely
from the experiments of Zhang~\cite{1994.Zhang} and Lindstr\"{o}m~\cite{1998.2.Lindstrom}. Second, there is some
literature (e.g.,~\cite{1995.Erdos.Sarkozy.Sos} and~\cite{1996.Ruzsa}) concerning the possible size of the
largest gap in a maximal Sidon set contained in $\{1,\dots,n\}$.  In short, we include this argument because
there is some reason to believe that this is a significant source of the error term in at least one case, and
because there is some reason to believe that improvement is possible.

\proofof{(vi)} The set
    \begin{equation*}
    S:= \Big[ 0, \floor{\frac g3} \Big) \cup \bigg\{ g - \floor{\frac g3} + 2 \Big[ 0, \floor{\frac g6} \Big)
    \bigg\} \cup \Big[ g, g+\floor{\frac g3} \Big) \cup \Big( 2g-\floor{\frac g3}, 3g-\floor{\frac g3} \Big]
    \end{equation*}
has cardinality $g+2\floor{g/3}+\floor{g/6}$, is contained in $[0,3g-\floor{g/3}]$, and has
    \[
    \G{S}{2} = g+2\floor{\frac g3} + \floor{\frac g6}.
    \]

We remark that this family of examples was motivated by the finite sequence
    \[
    S=(1,0,\tfrac12,1,0,1,1,1),\]
which has the property that its autocorrelations
    \[
    S\ast S = (1,0,1,2,\tfrac 14, 3, 3, 3, 3, 3, 3, 2, 3,2,1)
    \]
are small relative to the sum of its entries. In other words, the ratio of the $\ell^\infty$-norm of $S*S$ to
the $\ell^1$-norm of $S$ itself is small. If we could find a finite sequence of rational numbers for which the
corresponding ratio were smaller, it could possibly be converted into a family of examples that would improve
the lower bound for $\underline\rho(2g)$ in Theorem~\ref{R.Lower.Bounds.thm} for large $g$.

\subsection{Theorem~\ref{ExplicitValues} and Theorem~\ref{R.Lower.Bounds.thm}}

\begin{table}\footnotesize\begin{center}
 $
\begin{array}{ccccr}
   g    &  x & R(g,x) &                  \text{Witness}                   & R(g,x)/\sqrt{gx}             \\
\hline
   2    &  7 &   4    &                    \{1,2,5,7\}                    & \sqrt{8/7}     \approx 1.069    \\
   3    &  5 &   4    &                    \{1,2,3,5\}                    & \sqrt{16/15}   \approx 1.033    \\
   4    & 31 &  12    &         \{1,2,4,10,11,12,14,19,25,26,30,31\}      & \sqrt{36/31}   \approx 1.078    \\
   5    &  9 &   7    &               \{1,2,3,4,5,7,9\}                   & \sqrt{49/45}   \approx 1.043    \\
   6    & 20 &  12    &         \{1,2,3,4,5,6,9,10,13,15,19,20\}          & \sqrt{6/5}     \approx 1.095    \\
   7    & 15 &  11    &       \{1, 2, 3, 7, 8, 9, 10, 11, 12, 13, 15\}    & \sqrt{121/105} \approx 1.073    \\
   8    & 30 &  17    &  \{1,2,5,7,8,9,11,12,13,14,16,18,26,27,28,29,30\} & \sqrt{289/240} \approx 1.097    \\
   9    & 24 &  16    &     \{1,2,3,4,5,6,7,8,9,13,14,15,17,22,23,24\}    & \sqrt{32/27}   \approx 1.089    \\
   10   & 33 &  20 &\{1,2,4,5,6,7,8,9,10,11,13,15,20,21,22,23,30,31,32,33\}& \sqrt{40/33}  \approx 1.101    \\
   11   & 25 &  18 &\{ 1,2,3,4,5,11,12,13,14,15,16,17,18,19,20,21,23,25\} & \sqrt{324/275} \approx 1.085 \\
\end{array}
$ \caption{Important values of $R(g,x)$ and witnesses}\label{Witness.Table}
\end{center}\end{table}

Our plan is to employ the inequality of Theorem~\ref{C.constructions}(v) when $y$ is large, $f=2$, and $x\approx
\frac83 g$. In other words, we need nontrivial lower bounds for $C(2,y)$ for $y\to\infty$ and for $R(g,x)$ for
values of $x$ that are not much larger than $g$. The first need is filled by Theorem~\ref{C.constructions}(i),
(ii) or (iii), while the second need is filled by Theorem~\ref{C.constructions}(vi).

For any positive integers $x$ and $m\le\sqrt{n/x}$, the monotonicity of $R$ in the second variable gives
$R(2g,n) \ge R(2g,x(m^2-1)) \ge R(g,x)C(2,m^2-1)$ by Theorem~\ref{C.constructions}(v). If we choose $m$ to be
the largest prime not exceeding $\sqrt{n/x}$ (so that $m\gtrsim\sqrt{n/x}$ by the Prime Number Theorem), then
Theorem~\ref{C.constructions}(ii) gives $R(2g,n) \ge R(g,x)\cdot m \gtrsim R(g,x)\sqrt{n/x}$ for any fixed
positive integer $g$, and hence
    \[
    \sigma(2g) = \liminf_{n\to\infty} \frac{R(2g,n)}{\sqrt{gn}}
        \ge \liminf_{n\to\infty} \frac{R(g,x)\sqrt{n/x}}{\sqrt{gn}}
        =  \frac{R(g,x)}{\sqrt{gx}}.
    \]

The problem now is to choose $x$ so as to make $R(g,x)/\sqrt{gx}$ as large as we can manage for each~$g$. For
$g=2,3,\dots,11$, we use Table~\ref{R.Table} to choose $x=7$, 5, 31, 9, 20, 15, 30, 24, 33, and 25, respectively
(see Table~\ref{Witness.Table} for witnesses to the values claimed for $R(g,x)$). This yields
Theorem~\ref{ExplicitValues}.

We note that Habsieger \& Plagne~\cite{Habsieger.Plagne} have proven that $R(2,x)/\sqrt{2x}$ is actually
maximized at $x=7$. For $g>2$, we have chosen $x$ based solely on the computations reported in
Table~\ref{R.Table}. For general $g$, it appears that $R(g,x)/\sqrt{gx}$ is actually maximized at a fairly small
value of $x$, suggesting that this construction suffers from ``edge effects'' and is not best possible.

The first assertion of Theorem~\ref{R.Lower.Bounds.thm} is the immediate consequence of the obvious
$R(2g+1,n)\ge R(g,n)$. To prove the lower bound on $\sigma(2g)$, we set $x=3g-\floor{g/3}+1$ and appeal to
Theorem~\ref{C.constructions}(vi).

We remark that the above proof gives the more refined result
\begin{equation*}
R(2g,n) \ge \frac{11}{8\sqrt3}\sqrt{2gn} \,\Big( 1 + O\Big( g^{-1} + \Big( \frac ng \Big)^{\!(\alpha-1)/2}
\,\Big) \Big)
\end{equation*}
as $\frac ng$ and $g$ both go to infinity, where $\alpha<1$ is any number such that for sufficiently large $y$,
there is always a prime between $y-y^\alpha$ and $y$. For instance, we can take $\alpha=0.525$
by~\cite{2001.Baker.Harman.Pintz}. This clarification implies the final assertion of the theorem for even $g$,
and the obvious inequality $R(2g+1,n)\ge R(2g,n)$ implies the final assertion for odd $g$ as well.

\section{Significant Open Problems}
\label{section.OpenProblems}

It seems highly likely that
    \[
    \lim_{n\to \infty} \frac{R(g,n)}{\sqrt{n}}
    \]
is well-defined for each $g$, but this is known only for $g=2$ and $g=3$. It also seems likely that
    \[
    \lim_{n\to\infty} \frac{R(2g,n)}{R(2g+1,n)} =1.
    \]
The evidence so far is consistent with the conjecture $\lim_{g\to\infty} \sigma(g)=\sqrt{2}$.

One truly outstanding problem is to construct sets $S\subseteq\ZZ$ with $\G{S}{2}=4$ that are not the union of
two Sidon sets. In fact, all known constructions of sets with $\G{S}{2}\le g$ are not native, but are built up
by combining Sidon sets. It seems doubtful that this type of construction can be asymptotically densest
possible. The asymptotic growth of $R(4,n)$, or even of $C(4,n)$, is a major target.

As a computational observation, the set $S=B_{\langle 1,0 \rangle} \cup B_{\langle 1,1 \rangle} \cup B_{\langle
1,2 \rangle}$, where
    $$B_{\langle k_1,k_2 \rangle}:= \big\{a' \in [q^3-1] \colon \theta^{a'}-k_2 \theta^2 -k_1\theta\in
    \field{q}\big\}$$
and $\theta$ generates the multiplicative group of $\field{q^3}$, has the property that
    $$S\ast S\ast S(k)=\big|\left\{(s_1,s_2,s_3)\colon s_i\in S, \sum s_i =k\right\}\big| \le 81,$$
even when the sums are considered modulo $q^3-1$. As such, it seems likely that the generalizations of Bose's
and Singer's constructions given in this paper generalize further to give sets whose $h$-fold sums repeat a
bounded number of times. Proving this, however, will require a more efficient handling of systems of equations
than is presented in the current paper.

We direct our readers to the survey and annotated bibliography~\cite{2004.OBryant} for the current status of
these and other open problems related to Sidon sets.

\section*{Acknowledgement}

The authors thank Heini Halberstam for thoughtful readings of this manuscript and numerous helpful suggestions.
The first author was supported in part by grants from the Natural Sciences and Engineering Research Council. The
second author was supported by an NSF-VIGRE fellowship and by NSF grant DMS-0202460.


\begin{bibdiv}
\begin{biblist}
\bib{2001.Baker.Harman.Pintz}{article}{
    author={Baker, R. C.},
    author={Harman, G.},
    author={Pintz, J.},
     title={The difference between consecutive primes. II},
   journal={Proc. London Math. Soc. (3)},
    volume={83},
      date={2001},
    number={3},
     pages={532\ndash 562},
      issn={0024-6115},
    review={\MathReview{1851081}},
}

\bib{1942.Bose}{article}{
    author={Bose, R. C.},
     title={An affine analogue of Singer's theorem},
   journal={J. Indian Math. Soc. (N.S.)},
    volume={6},
      date={1942},
     pages={1\ndash 15},
    review={\MathReview{0006735}},
}

\bib{1962.Bose.Chowla}{article}{
    author={Bose, R. C.},
    author={Chowla, S.},
     title={Theorems in the additive theory of numbers},
   journal={Comment. Math. Helv.},
    volume={37},
      date={1962/1963},
     pages={141\ndash 147},
    review={\MathReview{0144877}},
}

\bib{1944.Chowla}{article}{
    author={Chowla, S.},
     title={Solution of a problem of Erd\"os and Turan in
            additive-number-theory},
   journal={Proc. Lahore Philos. Soc.},
    volume={6},
      date={1944},
     pages={13\ndash 14},
    review={\MathReview{0014116}},
}

\bib{Cilleruelo.Ruzsa.Trujillo}{article}{
    author={Cilleruelo, Javier},
    author={Ruzsa, Imre Z.},
    author={Trujillo, Carlos},
     title={Upper and lower bounds for finite $B\sb h[g]$ sequences},
   journal={J. Number Theory},
    volume={97},
      date={2002},
    number={1},
     pages={26\ndash 34},
      issn={0022-314X},
    review={\MathReview{1939134}},
}

\bib{1944.Erdos}{article}{
    author={Erd\"{o}s, P.},
     title={On a problem of Sidon in additive number theory and on some
            related problems. Addendum},
   journal={J. London Math. Soc.},
    volume={19},
      date={1944},
     pages={208},
    review={\MathReview{0014111}},
}

\bib{1995.Erdos.Sarkozy.Sos}{article}{
    author={Erd{\H{o}}s, P.},
    author={S\'{a}rk\"{o}zy, A.},
    author={S\'{o}s, V. T.},
     title={On sum sets of Sidon sets. II},
   journal={Israel J. Math.},
    volume={90},
      date={1995},
    number={1-3},
     pages={221\ndash 233},
      issn={0021-2172},
    review={\MathReview{1336324}},
}

\bib{1941.Erdos.Turan}{article}{
    author={Erd\"{o}s, P.},
    author={Tur\'{a}n, P.},
     title={On a problem of Sidon in additive number theory, and on some
            related problems},
   journal={J. London Math. Soc.},
    volume={16},
      date={1941},
     pages={212\ndash 215},
    review={\MathReview{0006197}},
}

\bib{2001.Green}{article}{
    author={Green, Ben},
     title={The number of squares and $B\sb h[g]$ sets},
   journal={Acta Arith.},
    volume={100},
      date={2001},
    number={4},
     pages={365\ndash 390},
      issn={0065-1036},
    review={\MathReview{1862059}},
}

\bib{1994.Guy}{book}{
    author={Guy, Richard K.},
     title={Unsolved problems in number theory},
    series={Problem Books in Mathematics},
 publisher={Springer-Verlag},
     place={New York},
      date={1994},
     pages={xvi+285},
      isbn={0-387-94289-0},
    review={\MathReview{1299330}},
}

\bib{Habsieger.Plagne}{article}{
    author={Habsieger, Laurent},
    author={Plagne, Alain},
     title={Ensembles $B\sb 2[2]$: l'\'etau se resserre},
  language={French},
   journal={Integers},
    volume={2},
      date={2002},
     pages={Paper A2, 20 pp. (electronic)},
    review={\MathReview{1896147}},
}

\bib{Kolountzakis}{article}{
    author={Kolountzakis, Mihail N.},
     title={The density of $B\sb h[g]$ sequences and the minimum of dense
            cosine sums},
   journal={J. Number Theory},
    volume={56},
      date={1996},
    number={1},
     pages={4\ndash 11},
      issn={0022-314X},
    review={\MathReview{1370193}},
}

\bib{1998.2.Lindstrom}{article}{
    author={Lindstr\"{o}m, Bernt},
     title={Finding finite $B\sb 2$-sequences faster},
   journal={Math. Comp.},
    volume={67},
      date={1998},
    number={223},
     pages={1173\ndash 1178},
      issn={0025-5718},
    review={\MathReview{1484901}},
}

\bib{Martin.OBryant}{article}{
    author={Martin, Greg},
    author={O'Bryant, Kevin},
    title ={Upper bounds for generalized Sidon sets},
    note  ={in preparation},
}

\bib{2004.OBryant}{article}{
    author = {O'Bryant, Kevin},
    title  = {A Complete Annotated Bibliography of Work Related to Sidon Sequences},
    journal= {Electron. J. Combin.},
    volume = {Dynamic Survey 11},
    date   = {2004},
    pages  = {39 pp (electronic)},
}

\bib{1993.Ruzsa}{article}{
    author={Ruzsa, Imre Z.},
     title={Solving a linear equation in a set of integers. I},
   journal={Acta Arith.},
    volume={65},
      date={1993},
    number={3},
     pages={259\ndash 282},
      issn={0065-1036},
    review={\MathReview{1254961}},
}

\bib{1996.Ruzsa}{article}{
    author={Ruzsa, Imre Z.},
     title={Sumsets of Sidon sets},
   journal={Acta Arith.},
    volume={77},
      date={1996},
    number={4},
     pages={353\ndash 359},
      issn={0065-1036},
    review={\MathReview{1414515}},
}

\bib{1932.Sidon}{article}{
    author = {Sidon, S.},
    title  = {Ein {S}atz \"uber trigonometrische {P}olynome und seine {A}nwendung in der {T}heorie der
    {F}ourier-{R}eihen},
    date   = {1932},
    volume = {106},
    pages  = {536 \ndash 539},
    journal= {Math. Ann.},
}

\bib{1938.Singer}{article}{
    author={Singer, James},
     title={A theorem in finite projective geometry and some applications to
            number theory},
   journal={Trans. Amer. Math. Soc.},
    volume={43},
      date={1938},
    number={3},
     pages={377\ndash 385},
      issn={0002-9947},
    review={\MathReview{1501951}},
}

\bib{1994.Zhang}{article}{
    author={Zhang, Zhen Xiang},
     title={Finding finite $B\sb 2$-sequences with larger $m-a\sp {1/2}\sb
            m$},
   journal={Math. Comp.},
    volume={63},
      date={1994},
    number={207},
     pages={403\ndash 414},
      issn={0025-5718},
    review={\MathReview{1223235}},
}
\end{biblist}
\end{bibdiv}
\end{document}